\documentclass[11pt]{article}
\usepackage{amsmath}
\usepackage{amssymb}
\usepackage{eucal}

\begin{document}

\baselineskip 16pt

\title{On the problem of existence and conjugacy  of injectors of generalized $\pi$-soluble groups\thanks{Research  is supported by
a NNSF grant of China (Grant \# 11301227) and Natural Science Foundation of Jiangsu Province (grant \# BK20130119)}}

\author{Xia Yin, Nanying Yang\\
{\small School of Science, Jiangnan University, Wuxi 214122, P. R. China}\\
 {\small E-mail: yangny@jiangnan.edu.cn}\\ \\
{N.T. Vorob'ev}\\
{\small Department of Mathematics, Masherov Vitebsk State University,
Vitebsk 210038, Belarus}\\
 {\small E-mail: vorobyovnt@tut.by}}

\date{}
\maketitle

\begin{abstract}  In this paper, we prove the existence and conjugacy of injectors of a generalized $\pi$-soluble groups for the Hartley class defined by a invariable Hartley function, and give a description of the structure of the injectors.

\end{abstract}

\let\thefootnoteorig\thefootnote
\renewcommand{\thefootnote}{\empty}


\let\thefootnote\thefootnoteorig

\section{Introduction}

Throughout this paper, all groups are finite and $p$ is a prime. $G$ always denotes a group, $|G|$ is the order of $G$, $\sigma (G)$ is the set of all primes dividing $|G|$,  $\pi$ denotes a set of some primes. Let $\mathbb{P}$ be the set of all primes, then let $\pi'=\mathbb{P}\setminus \pi$. We use $G_{\pi}$ denote a Hall $\pi$-subgroup of $G$.

Recall that a class $\mathfrak{F}$ of groups is called a Fitting class if $\mathfrak{F}$ is closed under taking normal subgroups and products of normal $\frak{F}$-subgroups.
As usual, we denote by $\frak{E}, \frak{S}, \frak{N}$ the classes of all groups, all soluble groups, all nilpotent groups, respectively; $\frak{E}_{\pi}, \frak{S}_{\pi}, \frak{N}_{\pi}$  denote the classes of all $\pi$-groups, all soluble $\pi$-groups, all nilpotent $\pi$-groups, respectively; and $\frak{S}^{\pi}$ and $\frak{N}^{\pi}$ to denote the classes of all $\pi$-soluble groups and all $\pi$-nilpotent groups, respectively. It is well known that all the above classes are Fitting classes.

From the definition of Fitting class, we see that for each non-empty Fitting class $\frak{F}$, every group $G$ has a unique maximal
normal $\frak{F}$-subgroup, which is called the $\frak{F}$-radical of $G$ and denoted by $G_{\frak{F}}$.

For any non-empty class $\frak{F}$ of groups, a subgroup $V$ of $G$ is said to be $\frak{F}$-maximal if $V\in \frak{F}$ and $U=V$ whenever $V\leq U\leq G$ and $U\in\frak{F}$.
A subgroup $V$ of a group $G$ is said to be an $\frak{F}$-injector of $G$ if $V\cap K$ is an
$\frak{F}$-maximal subgroup of $K$ for every subnormal subgroup $K$ of $G$. Clearly, every $\frak{F}$-injector of $G$ is an $\frak{F}$-maximal subgroup of $G$.

Fitting classes play an important role in the theory of groups. The importance of the theory of Fitting classes can firstly be seen in the following theorem proved by Fischer, Gasch$\ddot{u}$tz and Harley [1], which is in fact a generalization of the classical Sylow theorem and Hall theorem.

{\bf Theorem 1.1} (see \cite{1} or \cite[Theorem VIII,2.8]{DH}. Let $\frak{F}$ be a non-empty Fitting class. Then a soluble group possesses exactly one conjugacy class of $\frak{F}$-injectors.

 For a Fitting class $\frak{F}$, we let $\pi(\frak{F})=\bigcup _{G\in \frak{F}}\sigma(G)$.  Note that if $\frak{F}=\frak{N}_p$ is the Fitting class of all $p$-groups, then the $\frak{F}$-injectors of a group $G$ are Sylow $p$-subgroups of $G$; If $\frak{F}=\frak{E}_{\pi}$  and $G$ has a Hall $\pi$-soluble group, then the $\frak{F}$-injectors of $G$ are Hall $\pi$-subgroups of $G$ (see \cite[p.68, Ex.1]{Guo1} or \cite[p.238]{BalE}).

About the existence of $\frak{X}$-injector of $G$,  Shemetkov posed the following problem.

{\bf Problem 1.2} (Shemetkov \cite[Problem 11.117]{KT}). Let $\frak{X}$ be a Fitting class of soluble groups. Is it true that every finite non-soluble group possesses  an $\frak{X}$-injector ?

This problem has been resolved in \cite{Fors, Iranzo} for the Fitting classes $\frak{X}\in \{\frak{S}, \frak{S}_{\pi}, \frak{N}\}$.

In connection with this, a interesting problem is:  to find the conjugate class of injectors in any $\pi$-soluble groups. The first result in this direction is the following famous Chuchin's theorem \cite{Cun3}: A $\pi$-soluble $G$ possesses a Hall $\pi$-subgroup (that is, $\frak{E}_{\pi}$-injectors)  and any two Hall $\pi$-subgroups are conjugate.

As a development of Chuchin's theorem, Shemetkov and Guo proved the following

{\bf Theorem 1.3} (\cite[Theorem 2.2]{Shem1} and \cite{Guo2}).   Let $\frak{F}$ be a Fitting class, $\pi=\pi (\frak{F})$. If $G/G_{\frak{F}}$ is $\pi$-soluble, then $G$ has an $\frak{F}$-injector and any two $\frak{F}$-injector of $G$ are conjugate in $G$.

The product $\frak{FH}$ of two Fitting classes $\frak{F}$ and
$\frak{H}$ is the class $(G \mid G/G_{\frak{F}}\in \frak{H}).$ It is
well known that the product of any two Fitting classes is also a
Fitting class and the multiplication of Fitting classes satisfies
associative law (see [2, Theorem IX, 1.12(a)(c)]).

Following \cite{Hartl, Vorob}, a function $f$ : $\mathbb{P} \longrightarrow \{
\rm{nonempty \ Fitting \ classes}\}$ is called a Hartley function (or in brevity, H-function). A Fitting class $\mathfrak{F}$ is local if $$
\mathfrak{F} = \mathfrak{E}_{\pi (\mathfrak{F})}\bigcap (\bigcap_{p\in\pi
(\mathfrak{F})} f(p)\mathfrak{N}_p\mathfrak{E}_{p'})$$ for some
H-function $f$.

For a H-function $h$, let $\pi=Supp(h):=\{p\in \mathbb{P} : h(p)\neq \emptyset \}$, which is called the
support of the H-function $h$, and $LH(h)=\bigcap _{p\in \pi} h(p)\frak{E}_{p'}\frak{N}_p.$  Then, a Fitting class $\frak{F}$ is called Hartley class
if there exists a H-function $h$ such that $\frak{H}=LH(h)$.  In this case, $\frak{H}$ is said to be defined by the H-function $h$ or $h$ is an H-function of $\frak{H}$.

It is clear that $\frak{N} \subseteq LH(h)$. Hence, if $\frak{F}=LH(h),$ then  $\pi(\frak{F})=\mathbb{P}.$ Any Hartley class is a local Fitting class (see \cite[p.31]{Hauck2}). But the converse is not true in general (see \cite[p.207, 4.2]{Hartl}).
For two class functions $f$ and $h$, if $f(p)\subseteq h(p)$ for all $p\in \pi$, then we write that
$f\leq h$.

Concerning Fitting classes and injectors, the authors in \cite{GuoV} (see also \cite{Hartl}) posed the problem (in the universe  $\frak{S}$): {\sl Let $\mathfrak{F}$ be a local Fitting class of soluble groups, could we describe the $\mathfrak{F}$-injectors of a soluble group? } In view of Problem 1.2 and Theorem 1.3, the following more general question (in the class $\frak{E}$) naturally arise:

{\bf Problem 1.4.} For a local Fitting class $\frak{F}$ and a non-soluble group $G$ (in particular, a $\pi$-soluble group $G$), whether $G$ possesses an $\frak{F}$-injector and any two $\frak{F}$-injectors are conjugate ?

Note that there exist non-soluble groups $G$ and non-local Fitting classes $\frak{F}$ such that $G$ has no $\frak{F}$-injector (see, for example, \cite[7.1.3-7.1.4]{BalE}).

In this paper, by developing local method offered by Hartley \cite{Hartl}, we will resolve Problem 1.4 for partial $\pi$-soluble groups $G$ (in particular, for $\pi$-soluble groups) and the Hartley class defined by a invariable Hartley function. In fact, we will prove the following

{\bf Theorem 1.5. }  Support that $\frak{H}=LH(h)$ be a Hartley class defined by a invariable Hartley function $h$, that is, $h(p)=\frak{X}$ for all prime $p\in \pi=Supp(h),$ where $\frak{X}$ is some non-empty Fitting class, and $G\in \frak{X}\frak{S}^{\pi}$ (in particular, $G$ is a $\pi$-soluble group). Then the following statements hold:

(1) $G$ possesses an $\frak{H}$-injector and any two $\frak{H}$-injectors are conjugate in $G$;

(2) Every $\frak{H}$-injector $V$ of $G$ is of type $G_{\frak{X}\frak{E}_{\pi '}}L,$ where $L$ is the subgroup of $G$ such that $L/G_{\frak{X}}$ is an $\frak{N}_{\pi}$-injector of
some Hall $\pi$-subgroup of $G/G_{\frak{X}}.$

From Theorem 1.5, a series of new classical conjugate classes of injectors in any $\pi$-soluble group are obtained and the structure of injectors of some groups are described. For example, the following results directly follow from Theorem 1.5.

{\bf Corollary 1.5.1.}  Every $\pi$-soluble group possesses exactly one conjugacy class of $\frak{X}\frak{N}^{\pi}$-injectors.

{\bf Corollary 1.5.2.} Every $\pi$-soluble group has a injector with limited $\pi$-nilpotent length (that is, $(\frak{N}^{\pi})^k$-injector for any natural number $k$) and any two of them are conjugate.

If $\pi=\{ p\}$, then from Theorem 1.5 we have the following

{\bf Corollary 1.5.3.}  If $\frak{H}$ is an Hartley class of type $\frak{X}\frak{E}_{p'}\frak{N}_p$, where $\frak{X}$ is a nonempty Fitting class, and $G$ a group such that $G/G_{\frak{X}}$ is $p$-soluble, then $G$ possesses an $\frak{H}$-injector and any two of them are conjugate in $G$.

In the case when $\frak{X}=(1)$, where $(1)$ is the class of all identity groups, taking account of Theorem of Iranzo-Toress in \cite{IranzT}, we have

{\bf Corollary 1.5.4.} Every $p$-soluble group $G$ possesses exactly one conjugacy class of $p$-nilpotent injectors and each $p$-nilpotent injector is  a maximal $p$-nilpotent subgroup of $G$ containing the $p$-nilpotent radical of $G$.

Note that even if $\pi\subsetneqq \sigma (G)$ an $\frak{X}\subseteq \frak{S}$, the all statements of Theorem 1.5 and the corollaries are still new results for any soluble group $G$.

All unexplained notations and terminologies are standard. The reader is referred to \cite{DH, BalE, Guob} if necessary.

\section{Preliminaries}

Recall that $G$ is said to be $\pi$-soulbe if there exists a series of subnormal subgroups
$$1=G_0\lhd G_1\lhd G_2 \cdots \lhd G_m=G$$
such that every factor $G_i/G_{i-1}$ is either a $p$-group for some $p\in \pi$ or a $\pi '$-group, for $i=1, 2, \cdots , m$.
In particular, if $\pi=p$, then a $\pi$-soluble group is said to be $p$-soluble.

If  $\sigma (G)\subseteq \pi$, then $G$ is said to be a $\pi$-group. A subgroup $H$ of $G$ is said to be a Hall $\pi$-subgroup of $G$ if $\sigma (H)\subseteq \pi$ and $\sigma (|G:H|)\subseteq \pi '$.
$G$ is said to be $p$-nilpotent if $G$ has a normal Hall $p'$-subgroup. $G$ is said to be $\pi$-nilpotent if $G$ is $p$-nilpotent for all $p\in \pi$.

If $C_G(G_{\frak{F}})\subseteq G_{\frak{F}},$ then $G$ is said to be
$\frak{F}$-constrained. Note that if $\frak{F}=\frak{N}$ ($\frak{F}=\frak{N}^{\pi}$, resp.), then the $\frak{F}$-radical of $G$ is the
Fitting subgroup of $G$ (the $\pi$-Fitting subgroup, resp.), which is called nilpotent radical ($\pi$-nilpotent radical, resp.) and denoted by $G_{\frak{N}}$ or $F(G)$ ($G_{\frak{N}^{\pi}}$ or $F_{\pi}(G)$, resp.).
The maximal normal $\pi $-subgroup (the maximal normal $\pi '$-subgroup) of $G$ is said to be $\pi $-radical of $G$ and denoted by
$G_{\frak{E}_{\pi}}$ or $O_{\pi}(G)$ ($\pi '$-radical of $G$, and denoted by $G_{\frak{E}_{\pi '}}$ or $O_{\pi '}(G)$, resp.).

{\bf Lemma 2.1} (see \cite[Theorems 1.8.18 and 1.8.19]{Guo1} or \cite[Corollary 4.1.2]{Shem}). Support that $G\in \frak{S}^{\pi}.$ Then $G$ is $\frak{N}^{\pi}$-constrained, that is, $C_G(F_{\pi}(G))\leq F_{\pi}(G).$
In particular, if $G_{\frak{E}_{\pi '}}=1,$ then $G$ is $\pi$-constrained, that is, $C_G(O_{\pi}(G)\leq O_{\pi}(G).$

The following results is well known (see, for example, \cite[IX, Remarks (1.3)]{DH}).

{\bf Lemma 2.2.} Let $\frak{F}$ be a non-empty class of groups.

1) If $V$ is an $\frak{F}$-injector  of $G$ and $K\lhd G$, then $V\cap K$ is an  $\frak{F}$-injector of $K$;

2) If $V$ is an $\frak{F}$-injector  of $G$ and $\alpha : G \rightarrow \bar{G}$ is a isomorphism, then $\alpha (V)$ is an $\frak{F}$-injector  of $\bar{G}$;

3) If $V$ is $\frak{F}$-maximal subgroup of $G$ and $V\cap M$ is an $\frak{F}$-injector  of $M$ for any maximal normal subgroup $M$ of $G$, then $V$ is an $\frak{F}$-injector  of $G$;

4) If $V$ is an $\frak{F}$-injector  of $G$, then $G_{\frak{F}}\leq V$ and $V$ is an $\frak{F}$-maximal subgroup of $G$.

{\bf Lemma 2.3} \cite[Lemma IX, 1,1(a), Theorem IX, 1.12(b)]{DH}. Let $\frak{F}$ be a non-empty Fitting class. Then:

1) If $N$ is a subnormal subgroup of $G$, then $N_{\frak{F}}=N\cap G_{\frak{F}}$;

2) If $\frak{H}$ is a non-empty Fitting class, then the $\frak{H}$-radical of  $G/G_{\frak{F}}$ is $G_{\frak{FH}}/G_{\frak{F}}.$

{\bf Lemma 2.4} (see \cite{Cun3} or \cite[Chapter 5, Theorem 3.7]{Suz}). If $G\in \frak{S}^{\pi}$, then every $\pi$-subgroup of $G$ is contained some $\frak{E}_{\pi}$-injector of $G$ and any two $\frak{E}_{\pi}$-injector of $G$ are conjugate in $G$.

{\bf Definition 2.5.} Let $\pi =Supp(h)$, where $h$ is the support of some H-function $h$ of an Hartley class $\frak{H}$. Then $h$ is said to be

1)  integrated if $h(p)\subseteq {\frak H}$ for all $p\in \pi$;

2) full if $h(p)\subseteq h(q)\frak{E}_{p'}$ for all different primes $p, q\in \pi$;

3) full integrated if $h$ is full and integrated as well;

4) invariable if $f(p)=f(q)$ for all $p, q\in \pi$.

It is easy to see that every Hartley class can be defined by a integrated $H$-function, and a invariable H-function is full integrated (in fact, since $h(p)=h(q)$ for all $p, q\in \pi$, we have that $h(p)\subseteq h(p)\frak{E}_{q'}\subseteq h(q)\frak{E}_{q'}$, so $h(p)\subseteq \bigcap _ {q\in \pi}h(q)\frak{E}_{q'}\frak{N}_p=\frak{H}$).

\section{Proof of Theorem 1.5}

The proof of Theorem 1.5 consists of a large number of steps. The following 5 lemmas are the main steps of it.

{\bf Lemma 3.1.} Every $\pi$-soluble group $G$ possesses exactly one conjugacy class of $\pi$-nilpotent injectors, and each $\pi$-nilpotent injector of $G$ is the product of the $\pi '$-radical of $G$ and an $\frak{N}_{\pi}$-injector of some Hall $\pi$-subgroup of $G$.

{\bf Proof.}  We prove the lemma by induction on $|G|$. Let $M$ be any maximal normal subgroup of $G$. We consider the following two possible cases.

Case 1. The $\pi '$-radical $G_{\frak{E}_{\pi '}}$ of $G$ is a identity group, that is, $G_{\frak{E}_{\pi '}} = 1$.

Then $M_{\frak{E}_{\pi '}}=1$. By induction, $M$  possesses exactly one conjugacy class of $\pi$-nilpotent injectors  and every $\pi$-nilpotent injector of $M$ is an $\frak{N}_{\pi}$-injector of some Hall $\pi$-subgroup of $M$.

Let $F_1$ be an $\frak{N}_{\pi}$-injector of some Hall $\pi$-subgroup $M_{\pi}$ of $M$. Since every Hall $\pi$-subgroup $G_{\pi}$ of $G$ is soluble, $G_{\pi}$ has an  $\frak{N}_{\pi}$-injector $V$ and any two $\frak{N}_{\pi}$-injectors of $G_{\pi}$ are conjugate in $G_{\pi}$.
Since
$$M_{\pi}=M\cap G_{\pi}\lhd G_{\pi},$$
we have that $V\cap M_{\pi}$ is an $\frak{N}_{\pi}$-injector of $M$ by Lemma 2.2(1).  In view of the conjugacy of $\frak{N}_{\pi}$-injectors of $M$, we may assume that $F_1=V\cap M_{\pi}.$  Since a Hall $\pi$-subgroup of any $\pi$-nilpotent group is a nilpotent $\pi$-subgroup, every $\pi$-nilpotent injector of $G_{\pi}$ is an $\frak{N}_{\pi}$-injector of $G_{\pi}$. Hence, if we can prove that $V$ is a maximal $\pi$-nilpotent subgroup of $G$, then $V$ is a $\pi$-nilpotent injector of $G$ by Lemma 2.2(3), and so $V$ is an $\frak{N}_{\pi}$-injector of $G$.

Support that $V\leq V_1$ where $V_1$ is a maximal $\pi$-nilpotent subgroup of $G$. Since $G_{\frak{N}_{\pi}}$ and $(V_1)_{\frak{E}_{\pi '}}$ are normal in $V_1$, $[(V_1)_{\frak{G}_{\pi '}}, G_{\frak{N}_{\pi}}]\leq (V_1)_{\frak{E}_{\pi '}}\cap G_{\frak{N}_{\pi}}=1.$ Hence
$$(V_1)_{\frak{E}_{\pi '}}\leq C_G(G_{\frak{N}_{\pi}}).$$
Since $F_{\pi}(G)=G_{\frak{N}^{\pi}}=G_{\frak{N}_{\pi}}$ and $C_G(G_{\frak{N}^{\pi}})\leq G_{\frak{N}^{\pi}}$ by Lemma 2.1, we have that $(V_1)_{\frak{E}_{\pi '}}=1.$ This means that $V_1\in \frak{N}_{\pi}$ and so $V=V_1$ is a maximal $\pi$-nilpotent subgroup of $G$. This shows that the statement of the lemma holds in case 1.

Case 2.  $G_{\frak{E}_{\pi '}}\neq 1.$

Let $G_1=G/G_{\frak{E}_{\pi '}}.$  By Lemma 2.3(2), $(G_1)_{\frak{E}_{\pi '}}=G_{\frak{E}_{\pi '}}G_{\frak{E}_{\pi '}}/G_{\frak{E}_{\pi '}} =1.$  Hence by case 1, we have that $G_1$ possesses exactly one conjugate class of $\pi$-nilpotent injectors of type $(G_1)_{\frak{E}_{\pi '}}V_1,$ where $V_1$ is an $\frak{N}_{\pi}$-injector of some Hall $\pi$-subgroup of $G_1$, and the set of $\pi$-nilpotent injectors of $G_1$ are coincide with the set of $\frak{N}_{\pi}$-injectors of Hall $\pi$-subgroups $G_{\pi}G_{\frak{E}_{\pi '}}/G_{\frak{E}_{\pi '}}.$  But since $G_{\pi}$ is soluble, by Theorem 1.1 $G_{\pi}$ has $\frak{N}_{\pi}$-injector, V say. Then by Lemma 2.2(2), the subgroup $VG_{\frak{E}_{\pi '}}/G_{\frak{E}_{\pi '}}$ is an $\frak{N}_{\pi}$-injector of $G_{\pi}G_{\frak{E}_{\pi '}}/G_{\frak{E}_{\pi '}}$.  It follows that $VG_{\frak{E}_{\pi '}}$ is a $\pi$-nilpotent subgroup of $G$.

We now prove that $VG_{\frak{E}_{\pi '}}$ is a maximal $\pi$-nilpotent subgroup of $G$. Assume that $VG_{\frak{E}_{\pi '}}\leq F$ and $F$ is a maximal $\pi$-nilpotent subgroup of $G$. Then $F=F_{\frak{E}_{\pi '}}F_{\pi}$, where $F_{\pi}\in \frak{N}_{\pi}$ is a Hall $\pi$-subgroup of $F$. Without loss of generality, we may assume that $F_{\pi}\subseteq G_{\pi}$. Hence $V\leq F_{\pi}\leq G_{\pi}.$ But since the $\frak{N}_{\pi}$-injector $V$ is $\frak{N}_{\pi}$-maximal in $G_{\pi}$, we have that $V=F_{\pi}$. It follows from Lemma 2.3 and $G_{\frak{E}_{\pi '}}\lhd F$ that
$$(F/G_{\frak{E}_{\pi '}})_{\frak{E}_{\pi '}}=F_{\frak{E}_{\pi '}}/G_{\frak{E}_{\pi '}}.$$
Hence
$$(F/G_{\frak{E}_{\pi '}})/(F/G_{\frak{E}_{\pi '}})_{\frak{E}_{\pi '}}\simeq F/F_{\frak{E}_{\pi '}}\simeq F_{\pi}=V\in \frak{N}_{\pi}.$$
This shows that $F/G_{\frak{E}_{\pi '}}$ is $\pi$-nilpotent and $VG_{\frak{E}_{\pi '}}/G_{\frak{E}_{\pi '}}\leq F/G_{\frak{E}_{\pi '}}.$ Thus, $G_{\frak{E}_{\pi '}}V=F$ is a maximal $\pi$-nilpotent subgroup of $G$.

In order to prove that $G_{\frak{E}_{\pi '}}V$ is a $\pi$-nilpotent injector of $G$, by Lemma 2.2(3) we only need to prove that $G_{\frak{E}_{\pi '}}V\cap M$ is a $\pi$-nilpotent injector of $M$.

By induction, $M$ has a  $\pi$-nilpotent injector of type $M_{G_{\frak{E}_{\pi '}}}L$, where $L$ is an $\frak{N}_{\pi}$-injector of some Hall $\pi$-subgroup $M_{\pi}$ of $M$. Since
$$M_{\pi}=M\cap G_{\pi}\lhd G_{\pi},$$
and any two $\frak{N}_{\pi}$-injectors of $M_{\pi}$ are conjugate by Theorem 1.1, we may, without loss of generality, assume that $L=V\cap G_{\pi}.$
Since $G_{\frak{E}_{\pi '}}V\cap M\unlhd G_{\frak{E}_{\pi '}}V$ and $G_{\frak{E}_{\pi '}}V$ is $\pi$-nilpotent, $G_{\frak{E}_{\pi '}}V\cap M$ is $\pi$-nilpotent. But
$M_{\frak{E}_{\pi '}}L\leq  G_{\frak{E}_{\pi '}}V\cap M$ and the $\pi$-nilpotent injector $M_{\frak{E}_{\pi '}}L$ of $M$ is a maximal $\pi$-nilpotent subgroup of $M$. Therefore $M_{\frak{E}_{\pi '}}L=G_{\frak{E}_{\pi '}}V\cap M$. This shows that $G_{\frak{E}_{\pi '}}V\cap M$ is a $\pi$-nilpotent injector of $M$. Therefore, existence of $\pi$-nilpotent injector in a $\pi$-soluble group has been proved.

The conjugacy of $\pi$-nilpotent injectors follows from the cojugacy of $\frak{N}_{\pi}$-injectors of Hall $\pi$-subgroups of any $\pi$-soluble group. This shows that the lemma also holds in case 2.

The Lemma is proved.

In the case when $\pi=\{p\}$, by Lemma 3.1 we directly obtain the following

{\bf Corollary 3.2.} Every $p$-soluble group $G$ possesses exactly one conjugate class of $p$-nilpotent injectors, and every $p$-nilpotent injector of $G$ is of type $G_{\frak{E}_{p'}}P$, where $P$ is a Sylow $p$-subgroup of $G$ and $G_{\frak{E}_{p'}}$ is the $p'$-radical of $G$.

Support that $\frak{X}$ be some set of groups. We use ${\rm Fit}\frak{X}$ to denotes
the Fitting class generated by $\frak{X}$, that is, ${\rm
Fit}\frak{X}$ is the smallest Fitting class containing
$\frak{X}.$  For a class $\frak{F}$ of groups, we use $G^{\frak{F}}$ to denotes the $\frak{F}$-residual of  $G$.

{\bf Lemma 3.3.} Every Hartley class $\frak{H}$ can be defined  by
a full integrated $H$-function $h$, that is, $\frak{H}=LH(h)$ such that  $h(p)\subseteq
h(q)\frak{E}_{q'}\subseteq \frak{H},$ for all different primes $p, q\in \pi =Supp(h)$.

{\bf Proof.} Let $\frak{H}$ be a Hartley class. Then
$\frak{H}=LH(h_1)$, for some integrated $H$-function $h_1$. Following \cite{Hartl}, we
define:
$$\psi (p)=\{G \ |  \ G\simeq H^{\frak{E}_{p'}},
{\rm for \ some } \ H\in h_1(p)\},$$ for all $p\in \pi.$

Let $X$ be a group in $\psi (p).$ Then $X\simeq
Y^{\frak{E}_{p '}}$, for some group $Y\in h_1(p).$ Since
every Fitting class is closed with respect to normal subgroup,
$Y^{\frak{E}_{p '}}\in h_1(p)$ and so $X\in h_1(p).$
This shows that $\psi \leq h_1,$ and so
$$\psi (p)\frak{E}_{p'}\subseteq h_1(p)\frak{E}_{p'}.$$

If $Y_1\in h_1(p)\frak{E}_{p'}$, then
$Y_1/(Y_1)_{h_1(p)}\in \frak{E}_{p'}$ and $Y_1^{\frak{E}_{p'}}\leq (Y_1)_{h_1(p)}\in h_1(p).$ Moreover, since
$(Y_1^{\frak{E}_{p'}})^{\frak{E}_{p'}}=Y_1^{\frak{E}_{p '}},$ we have $Y_1^{\frak{E}_{p'}}\in \psi (p),$ that is, $Y_1\in \psi (p)\frak{E}_{p'}.$  Therefore, we obtain the following
equation:
$$\psi (p)\frak{E}_{p'}=h_1(p)\frak{E}_{p'}. \eqno(*)$$

Now, let $h$ be the $H$-function such that $h(p)={\rm Fit} (\psi
(p))$, for all $p\in \pi .$ Let
$$\frak{M}=\cap _{p\in \pi}h(p)\frak{E}_{p'}\frak{N}_{p}.$$

We now prove that $\frak{M}=\frak{H}.$ In fact, since $\psi \leq
h_1,$ we have $h\leq h_1$, and so $h(p)\frak{E}_{p'}\subseteq h_1(p)\frak{E}_{p'},$ for all $p\in \pi.$
Then, by the equation (*), we see that
$${\rm Fit} (h_1(p)\frak{E}_{p'})=h_1(p)\frak{E}_{p'}={\rm Fit} (\psi (p)\frak{E}_{p'}).$$ Therefore
$$h_1(p)\frak{E}_{p'}={\rm Fit} (\psi (p)\frak{E}_{p'})\subseteq ({\rm Fit} (\psi (p))\frak{E}_{p '}=h(p)\frak{E}_{p '}.$$
This shows that
$$h(p)\frak{E}_{p'}=h_1(p)\frak{E}_{p '},$$
for all $p\in \pi.$
Thus, we obtain that $\frak{M}=\frak{H}$. Moreover, since $h\leq h_1$ and $h_1$ is integrated $H$-function of $\frak{H}$, we have that
$h$ is a integrated $H$-function of $\frak{H}.$

Now, in order to complete the proof of the lemma, we only need to
prove that $h(p)\subseteq h(q)\frak{E}_{q'},$ for all different primes $p, q\in \pi.$
In fact, let $L\in h_1(p)$ and $p\neq q \in \pi.$
Clearly, $\frak{N}_{q}\subseteq
\frak{E}_{p '}$, so $L^{\frak{E}_{p'}}\leq
L^{\frak{N}_{q}}.$ However, since $L\in \frak{H}$, we have
$L/L_{h_1(q)\frak{E}_{q'}}\in \frak{N}_{q}.$
Hence $L^{\frak{N}_{q}}\leq  L_{h_1(q)\frak{E}_{q'}}.$  It follows that $L^{\frak{N}_{q}}\in h_1(q)\frak{E}_{q'}$. This shows that for every group in
$h_1(p)$, its $\frak{E}_{p'}$-residual is
contained in $h_1(q)\frak{E}_{q'}.$  Therefore, if $R\in
\psi (p),$ then $R\simeq V^{\frak{E}_{p'}}$ for some
group $V\in h_1(p)$ and thereby $R\in h_1(q)\frak{E}_{q'}.$ This induces that
$$\psi (p)\subseteq h_1(q)\frak{E}_{q'}.$$
Thus
$$h(p)={\rm Fit} (\psi (p))\subseteq {\rm Fit} (h_1(q)\frak{E}_{q'})=h_1(q)\frak{E}_{q'}=h(q)\frak{E}_{q'}.$$  This completes the proof of the lemma.

By Lemma 3.3, for any Hartley class $\frak{H}$, we may always
assume that $\frak{H}$ is defined by a full integrated $H$-function $h$. We
call the subgroup $G_h=\Pi _{p\in \pi}G_{h(p)}$ the {\sl $h$-radical} of $G$.

{\bf Lemma 3.4.} Let $\frak{H}=LH(h)$ for a full integrated $H$-function $h$, and $G$ a group such that $G/G_h$ is $\frak{N}^{\pi}$-constrained (in particular, $G/G_h$  is $\pi$-soluble). Then a subgroup $V$ containing $G_{\frak{H}}$ belongs to $\frak{H}$ if and only if $V/G_h$ is $\pi$-nilpotent.

{\bf Proof.} Assume that $V\in \frak{H}$ and $G_{\frak{H}}\subseteq V$. Then $V_{h(p)}\cap G_{\frak{H}}=(G_{\frak{H}})_{h(p)}=G_{h(p)},$ and so $[V_{h(p)}, G_{\frak{H}}]\subseteq G_{h(p)}.$
This implies that $V_{h(p)}\subseteq C_G(G_{\frak{H}}/G_{h(p)}),$ for all $p\in \pi$.

We first prove that $G_{\frak{H}}/G_h=F_{\pi}(G/G_h).$ Let $F_{\pi}(G/G_h)=L/G_h.$ Since $G_{\frak{H}}\in \frak{H}=\bigcap _{p\in \pi} h(p)\frak{E}_{p'}\frak{N}_p.$ and $(G_{\frak{H}})_{h(p)}=G_{h(p)},$  we have that $G_{\frak{H}}/G_{h(p)}$ is $p$-nilpotent for all $p\in \pi.$ Hence $G_{\frak{H}}/G_h$ is $\pi$-nilpotent. Consequently,  $G_{\frak{H}}/G_h\leq L/G_h$ and so $G_{\frak{H}}\leq L.$ On the other hand, since $L/G_h$ is $\pi$-nilpotent, by the isomorphism
$$L/L_{h(p)}G_h\simeq (L/G_h)/(L_{h(p)}G_h/G_h),$$
we have that $L/L_{h(p)}G_h \in \frak{E}_{p'}\frak{N}_p$ for all $p\in \pi.$ It follows that $(L/G_{h(p)}/(L_{h(p)}G_h/L_{h(p)})\in \frak{E}_{p'}\frak{N}_p$,
for all $p\in \pi.$

In order to prove that $L\leq G_{\frak{H}}$, we only need to prove that $L\in \frak{H}$. But it only need to prove that  $L_{h(p)}G_h/L_{h(p)}$ is a $p'$-group for all $p\in \pi$.
Since $G_h\unlhd L$, by Lemma 2.3,
$$G_{h(p)}=(G_h)_{h(p)}=G_h\cap L_{h(p)}\leq L_{h(p)}.$$
Let $q$ be an arbitrary prime in $\pi$ and $q\neq p$. Since
$$G_{h(p)}G_{h(q)}/G_{h(p)}\simeq G_{h(q)}/G_{h(q)}\cap G_{h(p)},$$
$G_{h(p)}G_{h(q)}/G_{h(p)}\simeq G_{h(q)}/(G_{h(q)})_{h(p)}.$
But since $h$ is full, we have that $h(q)\subseteq h(p)\frak{E}_{p'}$. Hence $G_{h(q)} \in h(p)\frak{E}_{p'},$ and so
$$G_{h(q)}/(G_{h(q)})_{h(p)} \in \frak{E}_{p'}.$$
This shows that $G_{h(p)}G_{h(q)}/G_{h(p)}$ is a $p'$-group for all different primes $p, q \in \pi$. Hence $G_h/G_{h(p)}\in \frak{E}_{p'}.$
Then by the isomorphism
$$L_{h(p)}G_h/L_{h(p)}\simeq G_h/G_h\cap L_{h(p)}\simeq (G_h/G_{h(p)})/(G_h\cap L_{h(p)}/G_{h(p)}),$$
we have that $L_{h(p)}G_h/L_{h(p)}$ is a $p'$-group for all $p\in \pi$. Thus $L\in \frak{H}$, and so $G_{\frak{H}}/G_h=F_{\pi}(G/G_h).$

Since $G/G_h$ is $\frak{N}^{\pi}$-constrained by the hypothesis and $G_{\frak{H}}/G_h=F_{\pi}(G/G_h)$, $C_{G/G_h}(G_{\frak{H}}/G_h)\leq G_{\frak{H}}/G_h$ and so $C_G(G_{\frak{H}}/G_h)\leq G_{\frak{H}}.$
But, clearly,  $C_G(G_{\frak{H}}/G_{h(p)})\leq C_G(G_{\frak{H}}/G_h).$ Hence $V_{h(p)}\subseteq G_{\frak{H}}$. Consequently, $V_{h(p)}=G_{h(p)}$ for all $p\in \pi$.
Then by $V\in \frak{H}=\bigcap _{p\in \pi} h(p)\frak{E}_{p'}\frak{N}_p$, we see that $V/G_{h(p)}=V/V_{h(p)}$ is $p$-nilpotent for all $p\in \pi.$ This implies that $V/G_h$ is $\pi$-nilpotent.

Conversely, if $V/G_h$ is $\pi$-nilpotent, then, with a similar argument (as the above proof of $L\leq G_{\frak{H}}$), we can see that $V\in \frak{H}.$
This completes the proof.

{\bf Lemma 3.5.} Let $\frak{H}=LH(h)$  for a full integrated $H$-function $h$, $\pi= Supp(h)\neq \emptyset$, and  $G$ be a group such that $G/G_h$ is $\frak{N}^{\pi}$-constrained and $\sigma(G_h)\subseteq \pi$.
If $V/G_h$ is a $\pi$-nilpotent injector of $G/G_h$, then $V$ is an $\frak{H}$-injector of $G$.

{\bf Proof.}  We prove the lemma by induction on $|G|$. Let $M$ be any maximal normal subgroup of $G$, and $p, q$ are different primes in $\pi$.

We first prove that $G_h/G_{h(q)}$ is a $q'$-group for all primes $q\in \pi.$
In fact, since $h(p)\subseteq h(q)\frak{E}_{q'}$  for all different primes $p, q\in \pi $, by the isomorphism
$$G_{h(q)}G_{h(p)}/G_{h(q)}\simeq G_{h(p)}/(G_{h(p)})_{G_{h(q)}},$$
we have that
$$G_{h(q)}G_{h(p)}/G_{h(q)}\in \frak{E}_{q'}.$$
This implies that $G_h/G_{h(q)}\in \frak{E}_{q'}$, that is, $G_h/G_{h(q)}$ is a $q'$-group.

Now let $M_h= \Pi _{p\in \pi}M_{h(p)}.$ Since
$$(G_h\cap M)G_{h(q)}/G_{h(q)}\simeq G_h\cap M/M_{h(q)},$$
$G_h\cap M/M_{h(q)}$ is a $q'$-group for all primes $q\in \pi.$ Hence
$$G_h\cap M/M_h\in \bigcap _{q\in \pi }\frak{E} _{q'} =\frak{E} _{\pi '}.$$
It follows from $\sigma (G_h)\subseteq \pi$ that
$$G_h\cap M/M_h\in \frak{E}_{\pi}\cap \frak{E} _{\pi '}=(1).$$
Thus $G_h\cap M=M_h.$

We consider the following two possible cases.

Case 1. $G_h\subseteq M$.

Then $G_h=M_h.$ Since $V/G_h$ is a $\frak{N}^{\pi}$-injector of $G/G_h$, $V\cap M/M_h$ is an $\frak{N}^{\pi}$-injector of $M/M_h$ by  Lemma 2.2(1).
Since $G/G_h$ is $\frak{N}^{\pi}$-constrained and the class of $\frak{N}^{\pi}$-constrained groups is a Fitting class by \cite[Theorem B(b)]{Izanzo}, $M/M_h$ is also $\frak{N}^{\pi}$-constrained.  Hence, by induction, $V\cap M$ is an $\frak{H}$-injector of $M$.

Assume that $V<V_1$ and $V_1$ is an $\frak{H}$-maximal subgroup of $G$. Since $V\cap M$ is an $\frak{H}$-maximal subgroup of $M$, we have that $V\cap M=V_1\cap M.$ Hence $V_1\cap M$ is an $\frak{H}$-injector of $M$ for any maximal subgroup $M$ of $G$. It follows from Lemma 2.2(3) that $V_1$ is a
$\frak{H}$-injector of $G$. But then $G_{\frak{H}}\leq V_1$, so $V_1/G_h$ is $\pi$-nilpotent by Lemma 3.4. As $V/G_h$ is an $\frak{N}^{\pi}$-injector of $G/G_h$, $V/G_h$ is a maximal $\pi$-nilpotent subgroup of $G/G_h$, which contradicts $V/G_h < V_1/G_h$. Thus $V=V_1$ is an $\frak{H}$-maximal subgroup of $G$. Hence by Lemma 2.2(3), $V$ is an $\frak{H}$-injector.

Case 2. $G_h\not\subseteq M$.

Then $G=G_hM$.  Since $V/G_h$ is a $\pi$-nilpotent injector of $G/G_h$ and
$$G/G_h\simeq M/G_h\cap M=M/M_h,$$
$V\cap M/M_h$ is an $\frak{N}^{\pi}$-injector of $M/M_h$ by Lemma 2.2(2). Then by induction, $V\cap M$ is an $\frak{H}$-injector of $M$. With a similar argument as in the case 1, we obtain that $V$ is an $\frak{H}$-injector of $G$.

The lemma is proved.

In view of the above Lemmas, we now may prove our main theorem.

{\bf Proof of Theorem 1.5.}  (1) Let $M$ be a maximal normal subgroup of $G$. Note that $h$ is a invariable H-function as $h(p)=\frak{X}$ for all $p\in \pi$. Hence $h$ is a full integrated H-function.
Since $G\in \frak{X}\frak{S}^{\pi}$, $G/G_{\frak{X}}$ is $\pi$-soluble. Therefore $G/G_{\frak{X}}$ is $\frak{N}^{\pi}$-constrained and $G/G_{\frak{X}}$ has a $\pi$-nilpotent injector $V/G_{\frak{X}}$ by Lemma 3.1. Clearly, $G_h=G_{\frak{X}}$, $M_h=M_{\frak{X}}$ and $G_h\cap M=M_h.$
Then, with the same arguments as in the proof of Lemma 3.5, we obtain that $V$ is an $\frak{H}$-injector of $G$.

Now we prove that if $V$ is an $\frak{H}$-injector of $G$, then $V/G_{\frak{X}}$ is a $\pi$-nilpotent injector of $G/G_{\frak{X}}$.

In fact, assume that $V$ is an $\frak{H}$-injector of $G$. Then  $V\cap S$ is an $\frak{H}$-maximal subgroup of $G$ for any subnormal subgroup $S$ of $G$.  In order to prove that $V/G_h$ is a $\pi$-nilpotent injector of $G/G_h$, we only need to prove that $V/G_{\frak{X}}\cap S/G_{\frak{X}}=(V\cap S)/G_{\frak{X}}$ is an $\frak{N}^{\pi}$-maximal subgroup of $G/G_{\frak{X}}$ for every subnormal subgroup $S/G_{\frak{X}}$ of $G/G_{\frak{X}}$.  Since $\frak{H}=LR(h)=\bigcap _{p\in \pi} \frak{X}\frak{E}_{p'}\frak{N}_p$ for all $p\in \pi$ and $V\cap S\in \frak{H}$, we have $(V\cap S)/G_{\frak{X}}\in \frak{N}^{\pi}$. Assume that $(V\cap S)/G_{\frak{X}}$ is not an $\frak{N}^{\pi}$-maximal subgroup of $G/G_{\frak{X}}$ and let $(V\cap M)/G_{\frak{X}} < D/G_{\frak{X}}$ and $D/G_{\frak{X}}$ is an $\frak{N}^{\pi}$-maximal subgroup of $G/G_{\frak{X}}$. Then clearly $D \in \bigcap _{p\in \pi} \frak{X}\frak{E}_{p'}\frak{N}_p=LR(h)=\frak{H}.$ But as $V\cap S$ is an $\frak{H}$-maximal subgroup of $G$, we have that $V\cap S=D.$ This contradiction shows that if $V$ is an $\frak{H}$-injector of $G$, then $V/G_{\frak{X}}$ is a $\pi$-nilpotent injector of $G/G_{\frak{X}}$.

Assume that $F$ is another $\frak{H}$-injector of $G$. Then $F/G_{\frak{X}}$ is a
$\pi$-nilpotent injector of $G/G_{\frak{X}}$ as above. Hence by Lemma 3.1, $V/G_{\frak{X}}$ and $F/G_{\frak{X}}$ are conjugate in $G/G_{\frak{X}}.$ This implies that $V$ and $F$ are conjugate in $G$. Hence we have (1).

(2) Let $V$ is a $\frak{H}$-injector of $G$. Then $V/G_{\frak{X}}$ is a $\pi$-nilpotent injector of $G/G_{\frak{X}}$. Hence by Lemma 3.1, we have that
$$V/G_{\frak{X}}=(G/G_{\frak{X}})_{\frak{E}_{\pi '}}(L/G_{\frak{X}}),$$
where $L/G_{\frak{X}}$ is an $\frak{N}_{\pi}$-injector of some Hall $\pi$-subgroup of $G/G_{\frak{X}}.$ But by Lemma 2.3,
$$(G/G_{\frak{X}})_{\frak{E}_{\pi '}}=G_{\frak{X}\frak{E}_{\pi '}}/G_{\frak{X}}.$$
This shows that $V=G_{\frak{X}\frak{E}_{\pi '}}L.$ Thus (2) holds.

The theorem is proved.

By Lemma 3.1, we known that every $\pi$-soluble group possesses exactly one conjugacy class of $\pi$-nilpotent injectors. Note also that every $\pi$-soluble group is $\frak{N}^{\pi}$-constrained. In connection with this, we put forward the following question:

{\bf Question 3.1.} Support that a group $G$ is $\frak{N}^{\pi}$-constrained. Is it true that $G$ possesses exactly one conjugacy class of $\pi$-nilpotent injectors?

\end{document}